\documentclass[12pt,a4paper]{amsart}
\usepackage[mathcal]{eucal}

\usepackage{amsmath, amssymb, xy}
\input xy
\xyoption{all}
\usepackage{pdflscape}
\usepackage{hyperref}
\usepackage[draft]{graphicx}

\numberwithin{equation}{section}

\newtheorem{thm}{Theorem}[section]

\newtheorem{lem}[thm]{Lemma}

\theoremstyle{definition}

\newtheorem*{remark*}{Remark}

\newcommand{\la}{{\lambda}}

\newcommand{\cA}{{\mathcal A}}

\newcommand{\End}{\text{End}}

\newcommand{\Hom}{\text{Hom}}
\newcommand{\rad}{\mathrm{rad}}

\title{The Cartan Matrix of a Centralizer Algebra}
\date{\today }

\author{Umesh V. Dubey}\address{The Institute of Mathematical Sciences, CIT campus, Taramani, Chennai 600113, India.}
\email{dubey@imsc.res.in}

\author{Amritanshu Prasad} \address{The Institute of Mathematical Sciences, CIT campus, Taramani, Chennai 600113, India.}
  \email{amri@imsc.res.in}

\author[Pooja Singla]{Pooja Singla${}^1$}\thanks{${}^1$Supported by the Center for Advanced Studies in Mathematics at Ben Gurion University} \address{Department of Mathematics, Ben Gurion
  University of the Negev, Beer-Sheva 84105 Israel}
  \email{pooja@math.bgu.ac.il}

\begin{document}
\maketitle
\begin{abstract}
The centralizer algebra of a matrix consists of those matrices that commute with it.
We investigate the basic representation-theoretic invariants of centralizer algebras, namely their radicals, projective indecomposable modules, injective indecomposable modules, simple modules and Cartan matrices.
With the help of our Cartan matrix calculations we determine their global dimensions.
Many of these algebras are of infinite global dimension. 
\end{abstract}

\section{Introduction}
\label{introduction} 
The centralizer algebra of a square matrix consists of all matrices that commute with it.
The centralizer algebra of a matrix shares many of the attributes of the matrix.
For example, a matrix is simple, semisimple or split if and only if its centralizer algebra has the respective properties. 

A fundamental invariant of a finite dimensional algebra over a field is its \emph{Cartan matrix}, introduced by \'Elie Cartan in \cite{MR1508194}.
It reflects, in a simple manner, some intrinsic properties of the category of finite dimensional modules for the algebra.
For example, if the algebra has finite global dimension, then its Cartan matrix has determinant $1$ or $-1$ \cite[Proposition~21]{MR0065544}.
The Cartan determinant conjecture says that  any finite dimensional algebra of finite global dimension has Cartan determinant  $1$.
For algebras with Loewy length less than equal to two, the finiteness of global dimension is reflected in the Cartan matrix, but for Loewy length 
three or more, knowledge of the Cartan matrix is not enough to determine finiteness of global dimension (see Burgess \emph{et al.} \cite{MR801315}). 
Another interesting feature of the Cartan determinant is its invariance under derived equivalence \cite[Proposition~1.5]{BokSko}.

In this article, the combinatorial data from the Jordan canonical form of a matrix is used to compute, in a direct and elementary manner, the Cartan matrix of its associated centralizer algebra.
A formula for its determinant follows easily.
This determinant is always a positive integer, and for most combinatorial data, it turns out to be greater than $1$.
It follows that the centralizer algebra has infinite global dimension.
When the determinant is $1$, the centralizer algebra turns out to be a sum of Auslander algebras, and therefore has global dimension $1$ or $2$.

\section{The Problem}
\label{problem}
We now move on to a precise formulation of the problem.
Definitions of the terms used here can be found in Section~\ref{prelims}.
Let $K$ be a perfect field and let $T$ be an $n\times n$-matrix with entries in $K$.
Let $\mathcal{A}$ denote the $K$-algebra of all matrices $B$ such that $TB=BT$. Then $\cA$ is called the centralizer algebra of $T$.
Let $P_1,\ldots,P_l$ be a complete set of representatives for the isomorphism classes of principal indecomposable $\cA$-modules (throughout this article, the term module should be understood to refer to a left module).
Each simple $\cA$-module is isomorphic to $D_i=P_i/\rad(P_i)$ for unique $i$.  
Given a finite dimensional $\cA$-module $M$ and a simple $\cA$-module $D$, let $[M:D]$ denote the number of composition factors in a composition series for $M$ that are isomorphic to $D$.
The \emph{Cartan matrix} of $\cA$ is the $l\times l$ matrix $C=(c_{ij})$ defined by
\begin{equation*}
  c_{ij}=[P_i:D_j].
\end{equation*}
The goal of this article is to compute the Cartan matrix of $\cA$.
On the way we also describe its radical and principal indecomposable modules.

\section{Preliminaries}
\label{prelims}
Let $A$ be a finite dimensional algebra over a field $K$ (not necessarily algebraically closed) with an identity. The radical of $A$, denoted $\rad(A)$, is the intersection of its
maximal left ideals, and turns out to be a two-sided ideal. The fact that $A$ is finite dimensional implies that $\rad(A)$ is nilpotent, that is, there exists 
a positive integer $m$ such that $(\rad(A))^m = (0)$. An algebra $A$ is called semisimple if its radical is zero. An example of a semisimple algebra is the algebra consisting of 
all matrices of order $n$ over $K$. Moreover, a direct sum of semisimple algebras is semisimple.  
The following characterization (see \cite[Section~24]{MR0144979}) of the radical of an algebra is useful.
\begin{lem}
\label{mod_nilpotent_semisimple} If $I$ is a two sided nilpotent ideal of $A$, then $I \subseteq \rad(A)$. If in addition, the algebra $A/I$ is semisimple then $I = \rad(A)$. 
\end{lem} 
An $A$-module $M$ is called {\it indecomposable} if for any $A$-modules $L$ and $N$, $M = L \oplus N$ implies either $L =0$ or $N = 0$. An $A$-module $M$ is called simple if it does
not have any proper submodules. A descending chain of submodules
\[
M = M_0 \supset M_1 \supset M_2 \cdots \supset M_t = (0)
\]
is called composition series if all the factor modules $M_{i}/M_{i+1}$ are simple $A$-modules. 
Clearly, any $A$-module $M$ with $\dim_{K}(M) < \infty$ has a composition series. 
\begin{thm}(Jordan holder Theorem)
Let $M$ be an $A$-module having two composition series
\[ 0 = M_0 \subset M_1 \subset \cdots \subset M_m = M
\]
\[
0 = N_0 \subset N_1 \subset \cdots \subset N_n = M 
\]
Then $m = n$ and there exists a permutation $\sigma$ of $\{0, 1, 2, \ldots, m \}$ such that for any $j \in \{0, 1, 2, \ldots, m-1\}$
\[
M_{j+1}/M_j \cong N_{\sigma(j+1)}/N_{\sigma(j)}. 
\] 
\end{thm}
For a proof, see \cite[Theorem~13.7]{MR0144979}.
Therefore the number of composition factors of $M$ isomorphic to a given simple $A$-module $S$ is a constant. We call this multiplicity of $S$ 
in $M$ and denote this by $[M:S]$.  

An element $e \in A$ is called idempotent if $e^2 = e$. Two idempotents $e_1$ and $e_2$ are called orthogonal if $e_1 e_2 = 0 = e_2 e_1$. An idempotent $e$ is called primitive
if it can not be written as direct sum of two non trivial orthogonal idempotents. Therefore an idempotent $e$ is primitive if and only if $Ae$ is an indecomposable $A$-module.
The algebra $A$, being finite dimensional over $K$, when considered as a left module over itself, decomposes into a finite direct sum of indecomposable modules, say $P_1$, $P_2$,..., $P_k$. 
These indecomposable $A$-modules are called the principal indecomposable $A$-modules (by \cite[Theorem~56.6]{MR0144979} these are the projective indecomposable $A$-modules as well). 
Further, for each principal indecomposable $P_i$ there exists a primitive idempotent $e_i$ such that $1 = \sum_{i=1}^k e_i$
and $e_i$'s are mutually orthogonal. By the Krull-Schmidt theorem~\cite[Theorem~14.5]{MR0144979}, the isomorphism classes of the principal indecomposable modules $P_1, P_2,\ldots, P_k$ are uniquely determined up to rearrangement. The radical of an $A$-module $M$, denoted $\rad(M)$, is the intersection of its maximal submodules. Suppose that $P_1, P_2, \ldots, P_l$ is a complete set of non-isomorphic principal indecomposable $A$-modules. The complete set of non-isomorphic simple and indecomposable injective $A$-modules can be determined from the principal indecomposable $A$-modules. 
\begin{lem} Let $P_1, P_2, \ldots, P_l$ be the complete set of representatives of isomorphism classes of
principal indecomposable $A$-modules and $e_1,e_2,\ldots,e_l$ be the corresponding primitive idempotents. 
\begin{enumerate}
\item Every simple $A$-module is isomorphic to one of the modules 
\[
S_1 = P_1/\rad(P_1), S_2 = P_2/\rad(P_2), \ldots, S_l = P_l/\rad(P_l).
\]
\item Every indecomposable injective $A$-module is isomorphic to one of the modules
\[
I_1 = \Hom_K(e_1 A, K), I_2 = \Hom_K(e_1 A, K), \ldots, I_l = \Hom_K(e_l A, K).
\]
\end{enumerate}
\end{lem}     
Any simple $A$-module has zero radical. So it is in fact $A/\rad(A)$-modules. Further the complete set of non-isomorphic simple $A$-modules coincides with the 
complete set of non-isomorphic simple $A/\rad(A)$-modules \cite[Theorem~25.24]{MR0144979}.
\section{Reduction to the primary case}
\label{sec:reduct-prim-case}
Observe that if $T$ and $T'$ are similar matrices, then their centralizer algebras are isomorphic. Therefore, for our purposes, we can always replace $T$ with a matrix similar to it. 
We will use $\bigoplus_i T_i$ to denote the block diagonal matrix whose diagonal blocks are $T_i$.
Then $T$ is similar to a matrix of the form 
\begin{equation*}
  \bigoplus_{p}T_p,
\end{equation*}
where $p$ ranges over the irreducible factors of the characteristic polynomial of $T$, and $T_p$ is a matrix whose characteristic polynomial is a power of $p$.
Moreover, $\cA$ has a decomposition into two-sided ideals
\begin{equation*}
  \cA=\bigoplus \cA_p,
\end{equation*}
where $\cA_p$ is the ring of matrices that commute with $T_p$.

Since $K$ is a perfect field, for each $p$ there exists a unique partition 
\[
\la = (\underbrace{\la_1, \ldots, \la_1}_{m_1 \mathrm{times}}, \underbrace{\la_2, \ldots, \la_2}_{m_2 \mathrm{times}}, \ldots, \underbrace{\la_r, \ldots, \la_r}_{m_r \mathrm{times}})
\]
with $\lambda_1<\lambda_2<\cdots <\lambda_r$, $m_1,\ldots,m_r$ positive integers, which we will abbreviate as
\begin{equation*}
  \lambda=(\lambda_1^{m_1},\ldots,\lambda_r^{m_r}),
\end{equation*}
such that
\begin{equation*}
  \cA_p\cong\End_{E[t]} \left(E[t]/(t^{\lambda_1})^{\oplus m_1}\oplus \cdots \oplus E[t]/(t^{\lambda_r})^{\oplus m_r}\right),
\end{equation*}
where $E$ is the algebraic field extension $K[t]/p(t)$ of $K$.
Note that the ring $\cA_p$ depends on $T_p$ only through $E$ and $\lambda$.
The partition $\lambda$ determines the shape of the Jordan canonical form of $\cA_p$.

We shall use the notation
\begin{equation*}
  E_\lambda = E[t]/(t^{\lambda_1})^{\oplus m_1}\oplus\cdots\oplus E[t]/(t^{\lambda_r})^{\oplus m_r},
\end{equation*}
which is an $E[t]$-module, and 
\begin{equation*}
  M_\lambda(E)=\End_{E[t]}(E_\lambda).
\end{equation*}
In particular $M_{1^m}(E)$ denotes the ring of $m\times m$ matrices over $E$.
\section{A block matrix representation of $M_\lambda(E)$}
An element $\vec{x}\in E_\lambda$ can be represented as a vector with entries $(\vec{x}_1,\ldots,\vec{x}_r)$, where $\vec{x}_j\in E[t]/(t^{\lambda_j})^{\oplus m_j}$.
Accordingly, if such a vector is represented as a column, an element $a\in M_\lambda(E)$ can likewise be represented by a matrix $a=(a_{ij})$, where
\begin{equation*}
  a_{ij}\in \Hom_{E[t]}\left(E[t]/(t^{\lambda_j})^{\oplus m_j}, E[t]/(t^{\lambda_i})^{\oplus m_i}\right).
\end{equation*}
Given $a=(a_{ij})$ and $b=(b_{ij})$ in $M_\lambda(E)$, the composition $ab$ has $(i,j)$ entry
\begin{equation*}
  a_{i1}b_{1j}+\cdots+a_{ir}b_{rj}.
\end{equation*}
Each of the above summands $a_{ik}b_{kj}$ is obtained by composition:
\begin{equation*}
  \xymatrix{
    E[t]/(t^{\lambda_j})^{\oplus m_j}\ar[r]^{b_{kj}} & E[t]/(t^{\lambda_k})^{\oplus m_k}\ar[r]^{a_{ik}} & E[t]/(t^{\lambda_i})^{\oplus m_i}
  }
\end{equation*}
\section{Computation of the radical}
\label{sec:computation-radical}
Consider $E_{\lambda_j}^{m_j}=E[t]/(t^{\lambda_j})^{\oplus{m_j}}$.
Reduction modulo $t$ gives a surjection $E_{\lambda_j}^{m_j}\to E^{m_j}$.
There is a corresponding reduction modulo $t$ for $M_{\lambda_j^{m_j}}(E)=\End_{E[t]}(E_{\lambda_j}^{m_j})$:
\begin{equation*}
  M_{\lambda_j^{m_j}}(E)\to M_{1^{m_j}}(E)
\end{equation*}
The kernel of above surjection, denoted $R_{\lambda_j}^{m_j}(E)$, consists of those endomorphisms of $E_{\lambda_j}^{m_j}$ whose image lies in $tE_{\lambda_j}^{m_j}$, hence 
is a two sided nilpotent
ideal of $M_{\lambda_j^{m_j}}(E)$. Therefore by Lemma~\ref{mod_nilpotent_semisimple}, it is the radical of
$M_{\lambda_j^{m_j}}(E)$. 
Let 
\begin{equation*}
  R_\lambda(E)=\{(a_{ij})\in M_\lambda(E) \mid a_{jj}\in R_{\lambda_j}^{m_j}(E)\text{ for } j=1,\ldots,r \}.
\end{equation*}
We claim that $R_\lambda(E)$ is the radical of $M_\lambda(E)$.
If $a\in R_\lambda(E)$ and $b\in M_\lambda(E)$, then the diagonal entries in the matrix of $ab$ are sums of terms of the form $a_{jk}b_{kj}$.
If $k>j$ then the image of $b_{kj}$ is contained in $tE_{\lambda_k}^{m_k}$.
If $k<j$ then the image of $a_{jk}$ is contained in $tE_{\lambda_j}^{m_j}$.
If $k=j$, then since $a\in R_\lambda(E)$, the image of $a_{jj}$ is contained in $tE_{\lambda_j}^{m_j}$.
Consequently, in all of these cases, the image of $a_{jk}b_{kj}$ is contained in $tE_{\lambda_j}^{m_j}$.
Therefore, $R_\lambda(E)$ is a right ideal.
The similar arguments shows that $R_\lambda(E)$ is a left ideal as well.
The quotient $\frac{M_\lambda(E)}{R_\lambda(E)}$ is a semisimple ring:
\begin{equation*}
  \frac{M_\lambda(E)}{R_\lambda(E)}=M_{1^{m_1}}(E)\oplus\cdots\oplus M_{1^{m_r}}(E).
\end{equation*}
It therefore remains only to show that $R_\lambda(E)$ is nilpotent.
For this, suppose that $a$ and $b$ are both in $R_\lambda(E)$.
Then the matrix entries of $ab$ are sums of terms of the form $a_{ik}b_{kj}$.
If $i\geq j$, then the sort of reasoning that was used earlier shows that $a_{ik}b_{kj}$ has image contained in $tE_{\lambda_i}^{m_i}$.
An inductive argument then shows that when $i\geq j$, then every product of $s$ elements in $R_\lambda(M)$ has $(i,j)$th entry whose image lies in $t^{s-1}E_{\lambda_i}^{m_i}$.
Therefore, the elements of $R_\lambda(M)^{\lambda_r}$ (recall $\la_r$ is largest part of partition) have all non-zero matrix entries strictly below the diagonal. 
A product of $r$ such elements is always zero.
Therefore $R_\lambda(E)^{\lambda_r r}=0$.
\section{Principal indecomposable and simple modules}
It follows from Sections~\ref{sec:computation-radical} and \ref{prelims} that there are $r$ isomorphism classes of simple $M_\lambda(E)$-modules.
These are represented by $D_1,\ldots,D_r$ where $D_j$, as an $E$-vector space is isomorphic to $E^{m_j}$, and $a\in M_\lambda(E)$ acts on it by the image of $a_{jj}$ in $M_{1^{m_j}}(E)$.
For each $1\leq i\leq r$ and $1\leq j \leq m_i$, let $e_{ij}$ denote the element of $M_\lambda(R)$ whose matrix has all entries $0$, except the $(i,i)$th entry, which as an element of $M_{\lambda_i^{m_i}}(E)$ has matrix with all entries zero and the $(j,j)$th entry equal to $1$.
Then
\begin{equation*}
  1=\sum_{i=1}^r\sum_{j=1}^{m_i} e_{ij},
\end{equation*}
and that the $e_{ij}$'s are pairwise orthogonal idempotents.
There is a matrix unit in $M_{\lambda_i^{m_i}}(E)$ which takes $e_{ij}$ to $e_{ij'}$ for any $1\leq j,j'\leq m_i$.
Therefore, all the modules $M_\lambda(E)e_{ij}$ are isomorphic for a fixed value of $i$.
Let $P_i=M_\lambda(E)e_{i1}$ for $i=1,\ldots,r$.
That $e_{i1}$ is a primitive idempotent follows from the corresponding fact for $M_{\lambda_i^{m_i}}$.
Since the $P_i/(R_\lambda(E)\cap P_i)$'s are pairwise non-isomorphic, so are the $P_i$'s.
It follows that $P_1,\ldots,P_r$ is a complete set of representatives for the isomorphism classes of principal indecomposable $M_\lambda(E)$-modules.
\section{The Cartan matrix}
The Cartan matrix of algebra $\cA$ is direct sum of Cartan matrices of the algebras $\cA_{p}$.
Therefore it now remains to calculate $[P_i :D_j]$ for all $1\leq i,j\leq r$.
By \cite[Theorem~54.16]{MR0144979}, 
\begin{equation*}
  [P_i:D_j]=\dim_E e_{j1}M_\lambda(E)e_{i1}.
\end{equation*}
For any $a\in M_\lambda(E)$, the $(m,n)$th  entry of $e_{j1}ae_{i1}$ is zero unless $m = j$ and $n = i$.
Think of $a_{ji}$ as an $m_i\times m_j$ matrix with entries in $\Hom_{E[t]}\left((E[t]/(t^{\lambda_i}),E[t]/(t^{\lambda_j})\right)$.
Then the $(j,i)$th entry of $e_{j1}ae_{i1}$ is $(a_{ji})_{11}$.
Moreover,
\begin{equation*}
  \dim_E \Hom_{E[t]}\left(E[t]/(t^{\lambda_i}),E[t]/(t^{\lambda_j})\right)=\min\{\lambda_i,\lambda_j\}
\end{equation*}
so that the Cartan matrix of $M_\lambda(E)$ is given by
\begin{equation*}
  c_{ij}=\min\{\lambda_i,\lambda_j\}, \text{ for all } 1\leq i,j\leq r.
\end{equation*}
In other words,
\begin{equation}
  \label{eq:1}
  C=
  \begin{pmatrix}
    \lambda_1 & \lambda_1 & \lambda_1 & \cdots & \lambda_1\\
    \lambda_1 & \lambda_2 & \lambda_2 & \cdots & \lambda_2\\
    \lambda_1 & \lambda_2 & \lambda_3 & \cdots & \lambda_3\\
    \vdots & \vdots & \vdots & \ddots & \vdots\\
    \lambda_1 & \lambda_2 & \lambda_3 & \cdots & \lambda_r
  \end{pmatrix}
\end{equation}
In particular $M_\lambda(E)$ has a single block.
Thus, $\cA$ has one block corresponding to each prime factor of the characteristic polynomial of $T$, and the Cartan matrix of each block is as described above.

\section{The Cartan Determinant}
\begin{lem}
  \label{lemma:det}
  The determinant of the matrix $C$ in (\ref{eq:1}) is given by
  \begin{equation*}
    \det C=\lambda_1(\lambda_2-\lambda_1)\cdots(\lambda_r-\lambda_{r-1})
  \end{equation*}.
\end{lem}
\begin{proof}
  Let $D(\lambda_1,\ldots,\lambda_r)$ denote the determinant of $C$.
  Subtracting the first row from each row of the matrix in (\ref{eq:1}) gives
  \begin{equation*}
    \begin{pmatrix}
      \lambda_1 & \lambda_1 & \lambda_1 & \cdots & \lambda_1\\
      0 & \lambda_2-\lambda_1 & \lambda_2-\lambda_1 & \cdots & \lambda_2-\lambda_1\\
      0 & \lambda_2-\lambda_1 & \lambda_3-\lambda_1 & \cdots & \cdots\\
      \vdots & \vdots & \vdots & \ddots & \vdots\\
      0 & \lambda_2-\lambda_1 & \lambda_3-\lambda_1 & \cdots & \lambda_r-\lambda_1
    \end{pmatrix}
  \end{equation*}
  whence
  \begin{equation*}
    D(\lambda_1,\ldots,\lambda_r)=\lambda_1 D(\lambda_2-\lambda_1,\ldots,\lambda_r-\lambda_1).
  \end{equation*}
  Applying the same process to the $(r-1)\times (r-1)$ determinant $D(\lambda_2-\lambda_1,\ldots,\lambda_r-\lambda_1)$ gives
  \begin{equation*}
    D(\lambda_2-\lambda_1,\ldots,\lambda_r-\lambda_1)=(\lambda_2-\lambda_1)D(\lambda_3-\lambda_2,\ldots,\lambda_r-\lambda_2),
  \end{equation*}
  and so on, until the required formula is obtained.
\end{proof}
Thus, $\det C$ is always a positive integer, and $\det C=1$ only when the $\lambda_i$'s are the consecutive integers $1,\ldots,r$.
Since finite global dimension implies that $\det C=\pm 1$, these are the only cases where $M_\lambda(E)$ has finite global dimension.
But these algebras are, by definition (see Section~\ref{sec:reduct-prim-case}) endomorphism algebras of additive generators for the category of $E[t]/(t^r)$-modules, and hence, by \cite[Proposition~5.2]{Ars}, have global dimension $2$, except for $r=1$, when the global dimension is $1$. We have
\begin{thm}
  The global dimension of $M_\lambda(E)$ is finite if and only if $\lambda$ is of the form $(1^{m_1},2^{m_2},\ldots,r^{m_r})$ for some positive integer $r$.
  When $r=1$, the global dimension is $1$.
  When $r\geq 2$, the global dimension is $2$.
\end{thm}

\end{document}